\documentclass[dvips]{article}
\usepackage{graphicx}
\catcode`\@=11
\newtheorem{lemma}{Lemma}
\newtheorem{proposition}{Proposition}

\newtheorem{example}{Example}

\newtheorem{corollary}{Corollary}
\def\demo{\noindent{\bf Proof .-}}
\def\section{\@startsection {section}{1}{\z@}{-3.5ex plus -1ex
minus-.2ex}{2.3ex plus .2ex}{\normalsize\bf}}

\begin{document}
\begin{center}
{\Large\bf \textsc{A note on the edge ideals of Ferrers graphs}}
\end{center}
\vskip.5truecm
\begin{center}
{Margherita Barile\footnote{Partially supported by the Italian Ministry of Education, University and Research.}\\ Dipartimento di Matematica, Universit\`{a} di Bari, Via E. Orabona 4,\\70125 Bari, Italy}\footnote{e-mail: barile@dm.uniba.it, Fax: 0039 080 596 3612}
\end{center}
\vskip1truecm
\noindent
{\bf Abstract} We determine the arithmetical rank  of every edge ideal of a Ferrers graph. 

\section*{Introduction} 
In a Noetherian commutative ring with identity, every ideal is finitely generated; in particular, for every ideal $I$, there are finitely many elements $f_1,\dots, f_s\in R$ such that Rad\,$I =$\,Rad\,$(f_1,\dots, f_s)$.  The least such number $s$ is called the {\it arithmetical rank} of $I$, denoted ara\,$I$. In this paper we determine the arithmetical rank of a special class of ideals generated by squarefree monomials in a polynomial ring over a field. These are the edge ideals of certain bipartite graphs, the so-called {\it Ferrers graphs}, and have been considered for the first time by Corso and Nagel. This paper is intended as an application of their joint work \cite{CN}, which contains an extensive treatment of the algebraic properties of the so-called {\it Ferrers ideals}. These admit a combinatorial description based on the {\it Ferrers diagrams} (also known as {\it Young diagrams}). We show that for every such ideal, the arithmetical rank  always equals its cohomological dimension (equivalently, its projective dimension). We also give a simple characterization of the cases where the arithmetical rank coincides with the height.

\section*{On the arithmetical rank of Ferrers ideals} 
Let $R$ be a commutative ring with identity.  We will use the following preliminary result, which is due to Schmitt and Vogel. 

\begin{lemma}\label{lemma}{\rm [\cite{SV}, p.249]} Let $P$ be a finite subset of elements of $R$. Let $P_1,\dots, P_r$ be subsets of $P$ such that
\begin{list}{}{}
\item[(i)] $\bigcup_{i=1}^rP_i=P$;
\item[(ii)] $P_1$ has exactly one element;
\item[(iii)] if $p$ and $p'$ are different elements of $P_i$ $(1<i\leq r)$ there is an integer $i'$ with $1\leq i'<i$ and an element in $P_{i'}$ which divides $pp'$.
\end{list}
\noindent
We set $q_i=\sum_{p\in P_i}p^{e(p)}$, where $e(p)\geq1$ are arbitrary integers. We will write $(P)$ for the ideal of $R$ generated by the elements of $P$.  Then we get
$${\rm Rad}\,(P)={\rm Rad}\,(q_1,\dots,q_r).$$
\end{lemma}
\noindent
We now introduce the class of ideals that we are going to examine. 
Given a positive integer $n$, and an $n$-uple $\lambda=(\lambda_1,\dots, \lambda_n)$ of positive integers such that $m=\lambda_1\geq\cdots\geq \lambda_n$, the {\it Ferrers graph} $G=G_{\lambda}$ associated with $\lambda$ is the bipartite graph on the vertex set $\{x_1,\dots, x_n, y_1,\dots, y_m\}$ whose set of edges is 
\begin{equation}\label{EG}E(G)=\{(x_i,y_j)\vert 1\leq i\leq n,1\leq j\leq\lambda_i\}.\end{equation}
\noindent
Suppose that $x_1,\dots, x_n, y_1,\dots, y_m$ are indeterminates over the field $K$. The {\it edge ideal} of $G$ in the polynomial ring $R=K[x_1,\dots, x_n, y_1,\dots, y_m]$ is the squarefree monomial ideal 
\begin{equation}\label{IG} I(G)=\big(\{x_iy_j\vert (x_i,y_j)\in E(G)\}).\end{equation}
\noindent
Corso and Nagel \cite{CN} call this a {\it Ferrers ideal}: their paper \cite{CN} is entirely devoted to the determination of the main algebraic invariants of this class of ideals. 
In \cite{CN}, Section 2, it is shown that  
\begin{equation}\label{redundant}I(G)=\bigcap_{i=1}^{n+1}(x_1,\dots, x_{i-1}, y_1, \dots, y_{\lambda_i}),\end{equation}
\noindent
where, by convenience of notation, we have set $\lambda_{n+1}=0$. 
Set $c_1=1$, and suppose that 
$$\lambda_1=\cdots=\lambda_{c_2-1}>\lambda_{c_2}=\cdots=\lambda_{c_3-1}>\lambda_{c_3}=\cdots=\lambda_{c_k-1}>\lambda_{c_k}=\cdots=\lambda_n.$$
Finally set $c_{k+1}=n+1$. Then a minimal prime decomposition of $I(G)$ can be obtained as follows, by omitting redundant terms from (\ref{redundant}):
\begin{equation}\label{irredundant}I(G)=\bigcap_{i=1}^{k+1}(x_1,\dots, x_{c_i-1}, y_1, \dots, y_{\lambda{c_i}}).\end{equation}
\noindent
The graph $G$ can be associated with the {\it Ferrers diagram} having $n$ rows and $m$ columns, and whose rows have lengths $\lambda_1,\dots, \lambda_n$. 
Note that $k$ is the number of outer corners of the Ferrers diagram associated with $G$.
\begin{proposition}\label{proposition} Let $\mu=\max_{j=1,\dots,n}\{\lambda_j+j-1\}$. For all $i=1,\dots,\mu$, set 
$$q_i=\sum_{\scriptstyle{r+s=i+1}\atop{x_ry_s\in I(G)}}x_ry_s.$$
\noindent Then $I(G)={\rm Rad}\,(q_1,\dots, q_{\mu})$. 
\end{proposition}
\demo Let $P$ be the set of minimal monomial generators of $I(G)$ given in (\ref{IG}), and, for all $i=1,\dots, \mu$, let $P_i=\{x_ry_s\vert r+s=i+1, x_ry_s\in P\}$. It suffices to show that $P$, and $P_1,\dots, P_{\mu}$ fulfill the assumption of Lemma \ref{lemma}. Condition (i) is true by definition of $I(G)$ (see (\ref{EG}) and (\ref{IG})). Since $P_1=\{x_1y_1\}$, condition (ii) is satisfied, too. We show (iii). 
Let $p=x_ry_s$ and $p'=x_uy_v$ be distinct elements of $P_i$ for some index $i=2,\dots, \mu$. Then $r+s=u+v=i+1$, and, consequently, $r\ne u$, $s\ne v$.  We may assume that $r<u$ and $s>v$. Then $v\leq \lambda_u\leq\lambda_r$, so that $p''=x_ry_v\in I(G)$. More precisely, $x_ry_v\in P_{r+v-1}$, where $r+v-1<r+s-1=i$. Since $p''$ divides $pp'$, this proves (iii) and completes the proof.
\par\medskip\noindent
Proposition \ref{proposition} shows that 
\begin{equation}\label{ara}{\rm ara}\,I(G)\leq \max_{j=1,\dots,n}\{\lambda_j+j-1\}.\end{equation}  
As is well known (see, e.g., \cite{H}, Example 2, pp.~414--415), a lower bound for the arithmetical rank of an ideal $I$  of $R$ is given by its {\it local cohomological dimension}
$${\rm cd}\, I=\max\{i\vert H_I^i(R)\ne 0\},$$
\noindent
where $H_I^i(R)$ denotes the $i$th local cohomology module of $R$ with respect to $I$. On the other hand, according to \cite{L}, Theorem 1, if $I$ is a squarefree monomial ideal, then cd\,$I=\,$pd\,$I$, where pd denotes the projective dimension, i.e., the length of any minimal free resolution of $I$ over $R$.
Thus we have that
\begin{equation}\label{cd}{\rm cd}\,I(G)={\rm pd}\,I(G)\leq\,{\rm ara}\, I(G).\end{equation}
 In \cite{CN} all Betti numbers of $I(G)$ are expressed by a closed formula in terms of $\lambda$. In particular we have that 
\begin{equation}\label{pd}{\rm pd}\,I(G)=\max_{j=1,\dots,n}\{\lambda_j+j-1\}.\end{equation}
\noindent
From (\ref{ara}), (\ref{cd}) and (\ref{pd}) we derive the following
\begin{corollary}\label{corollary} $${\rm ara}\,I(G)=\,{\rm cd}\,I(G)=\,{\rm pd}\, I(G)=\max_{j=1,\dots,n}\{\lambda_j+j-1\}.$$
\end{corollary}
 Recall that, for every ideal $I$ in $R$, ${\rm ht}\,I\leq{\rm ara}\,I$, where ht denotes the height. 
 If ${\rm ara}\,I={\rm ht}\,I$, $I$ is called a {\it set-theoretic complete intersection}. 
\begin{corollary}\label{corollary2} $I(G)$ is a set-theoretic complete intersection if and only if $\lambda=(m, m-1,\dots, 2,1)$.  
\end{corollary}
\demo It is well known that every set-theoretic complete intersection has pure height. Hence, by (\ref{irredundant}), we may assume that $m=n$. Moreover, in view of (\ref{redundant}),  ht\,$I(G)=\min_{j=1\,\dots, n+1}\{\lambda_j+j-1\}$. Note that $m=\lambda_1+1-1=\lambda_{n+1}+n+1-1=n$.  Therefore ht\,$I(G)=\min_{j=1\,\dots, n}\{\lambda_j+j-1\}$. On the other hand, from Corollary \ref{corollary} we have that ara\,$I(G)=\max_{j=1\,\dots, n}\{\lambda_j+j-1\}$. Hence  $I(G)$ is a set-theoretic complete intersection if and only if the minimum $\lambda_j+j-1$ coincides with the maximum $\lambda_j+j-1$, in other words, if and only if all $\lambda_j+j-1$ are equal. This occurs if and only if, for all indices $i=1,\dots, n$, it holds $m=\lambda_1=\lambda_i+i-1$, i.e., $\lambda_i=m-i+1$.  This proves the claim. 
\par\medskip\noindent
The {\it if} part of Corollary \ref{corollary2} had already been proven, by different methods, in \cite{B1}, Corollary 2. 
\begin{example}{\rm Consider the edge ideal of $G=G_{\lambda}$, where $\lambda=(6,4,4,2,1)$: 
\begin{eqnarray*}
&I(G)=(x_1y_1, x_1y_2, x_1y_3, x_1y_4, x_1y_5, x_1y_6, x_2y_1, x_2y_2, x_2y_3, x_2y_4,\hphantom{xxxxxxxx}\\
 &\hphantom{xxxxxxxxxxxxxxxxxxxxxxxxx}x_3y_1, x_3y_2, x_3y_3, x_3y_4, x_4y_1, x_4y_2, x_5y_1).
\end{eqnarray*}
\noindent
We have that ht\,$I(G)=5$ and ara\,$I(G)=$\,cd\,$I(G)=6$. Proposition \ref{proposition} yields:
\begin{eqnarray*}
&I(G)=\,{\rm Rad}\,(x_1y_1,\ x_1y_2+x_2y_1,\ x_1y_3+x_2y_2+x_3y_1,\hphantom{xxxxxxxxxxxxxx}\\ 
&x_1y_4+x_2y_3+x_3y_2+x_4y_1,\
x_1y_5+x_2y_4+x_3y_3+x_4y_2+x_5y_1,\ x_1y_6+x_3y_4).
\end{eqnarray*}
\noindent
Note that the sums are taken along the ascending diagonals of the associated Ferrers diagram.
In particular, the arithmetical rank can be described combinatorially as the number of such diagonals which are, at least partially, present in the diagram.} 
\end{example}

\end{document}